\documentclass[11pt]{article}
\usepackage{authblk}

\usepackage{authblk} 

\usepackage{amsmath}
\usepackage{amsfonts}
\usepackage{amssymb,bm}

\usepackage[normalem]{ulem} 

\usepackage{amsmath,xcolor,ulem}
\usepackage{amsfonts}
\usepackage{amssymb,todonotes}
\usepackage{amsthm}
\usepackage{graphicx}
\usepackage{marginnote}

\usepackage[colorlinks=true, allcolors=blue]{hyperref}

\usepackage[top=2.8cm,bottom=2.8cm,left=2.5cm,right=2.5cm]{geometry}
\usepackage{float}
\usepackage[algoruled]{algorithm2e}
\usepackage[font=footnotesize,width=\linewidth]{caption} 

\usepackage{graphicx}
\usepackage{amssymb}
\usepackage{amsthm}
\usepackage{amsmath}
\usepackage{algorithmic}

\usepackage{subcaption}
\usepackage{caption}
\usepackage{float}
\usepackage{geometry}
\usepackage{mathrsfs} 
\usepackage{diagbox}
\newtheorem{example}{Example}

\newtheorem{remark}{Remark}

\newcommand{\fitmat}[1]{{#1}}
\newcommand{\fitvec}[1]{{#1}}

\newcommand{\fMeps}{{M}_{\varepsilon}}
\newcommand{\fMkap}{{M}_{\kappa}}

\newcommand{\fMmu}{{M}_{\mu}}

\newcommand{\Gr}{G}
\newcommand{\C}{C}

\DeclareMathOperator{\diag}{diag}

\usepackage{xcolor}

\title{Structure-Preserving Coupling and Decoupling \\ of Port-Hamiltonian Systems}

\author{Matthias Ehrhardt\footnote{ \href{mailto:ehrhardt@uni-wuppertal.de}{ehrhardt@uni-wuppertal.de}} ,
Michael Günther\footnote{ \href{mailto:guenther@uni-wuppertal.de}{guenther@uni-wuppertal.de}} ,
Daniel \v{S}ev\v{c}ovi\v{c}\footnote{Corresponding author, \href{mailto:sevcovic@fmph.uniba.sk}{sevcovic@fmph.uniba.sk}}  
}

\affil{University of Wuppertal, Chair of Applied and Computational Mathematics,\\
Gaußstrasse 20, 42119 Wuppertal, Germany}

\affil{Comenius University Bratislava,
Department of Applied Mathematics and Statistics,\\
Mlynsk\'a dolina, 842 48 Bratislava, Slovakia}

\begin{document}
\maketitle

\begin{tikzpicture}[remember picture,overlay]
	\node[anchor=north east,inner sep=20pt] at (current page.north east)
	{\includegraphics[scale=0.2]{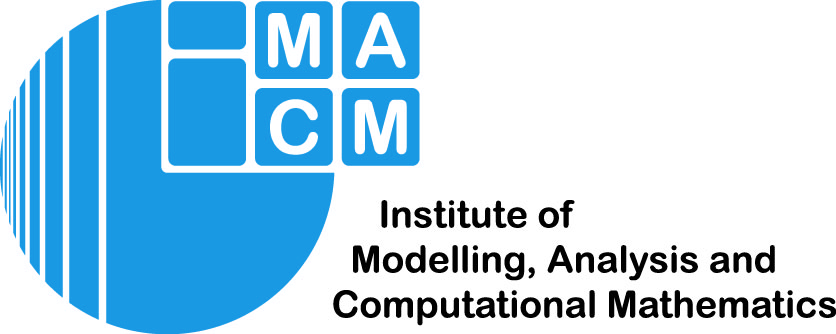}};
\end{tikzpicture}

\begin{abstract}
The port-Hamiltonian framework is a structure-preserving modeling approach that preserves key physical properties such as energy conservation and dissipation. When subsystems are modeled as port-Hamiltonian systems (pHS) with linearly related inputs and outputs, their interconnection remains port-Hamiltonian. This paper introduces a systematic method for transforming coupled port Hamiltonian ordinary differential equations systems (pHODE) into a single monolithic formulation, and for decomposing a monolithic system into weakly coupled subsystems. The monolithic representation ensures stability and structural integrity, whereas the decoupled form enables efficient distributed simulation via operator splitting or dynamic iteration.
\end{abstract}

\begin{minipage}{0.9\linewidth}

\medskip

\end{minipage}

\section{Introduction}\label{sec1}

In today’s rapidly evolving technological landscape, the ability to model and analyze complex systems is more crucial than ever. Port-Hamiltonian Systems (pHS) integrate Hamiltonian mechanics with network theory to provide unparalleled insight 
into multi-component systems. First, pHS exhibit a Hamiltonian structure, i.e., pHS provide an 
\textit{energy-based modeling approach}, essential for understanding the dynamics of modern systems (see related papers \cite{buchfink26}, \cite{Cervera2007}, \cite{Jaeschke22a}, \cite{Jaeschke22b} and book \cite{jacobbook}).
Next, they allow for \textit{dissipation} in the system. This accurately models energy losses, ensuring realistic simulations and designs. They have \textit{ports} that allow for input/output connections to the outside world, i.e. this approach facilitates energy exchange with the external environment. Finally, pHS have an \textit{interconnection property}: multiple subsystems can be coupled to a monolithic joint pH system, preserving the overall energy structure and increasing reliability of the system.

In this work, we investigate a) how coupled port-Hamiltonian ordinary differential equation (pHODE) systems can be formulated as a single joint pHODE and b) vice-versa, how a single joint pHODE system can be decoupled into subsystems of pHODE, and generalize some ideas mentioned  in \cite{Bartel23}. While a) allows for monolithic simulation with excellent stability properties, b) is the starting point for an efficient distributed simulation via operator splitting and/or dynamic
iteration.

The paper is organized as follows: after an introduction, in Section~\ref{sec:pHODE} we define the setting of port-Hamiltonian ODE systems. Section~\ref{sec:coupling} deals with the coupling of pHODE subsystems into a joint pHODE, and Section~\ref{section-decoupling1} discusses the decoupling of a single pHODE system into pHODE subsystems.
Section~\ref{sec:advantages} discusses briefly the advantages and disadvantages of both approaches. 
We conclude with final remarks in Section~\ref{sec:conclusion}. 

\section{Port-Hamiltonian ODE systems}\label{sec:pHODE}
Following the lines of~\cite[Definition 1]{MehrmannMorandin}, we  
define the general setting of port-Hamiltonian ordinary differential equations (pHODEs).
Consider a time interval $I \subset \mathbb{R}$, a state space $\mathcal{X} \subset \mathbb{R}^n$.
A pHODE is a (possibly implicit) system of ordinary differential equations (ODEs) of the form:
\begin{subequations}\label{pHODE}
    \begin{align}
        E(x) \dot{x} 
        & = \bigl( J(x)-R(x) \bigr)\, z(x) + \bigl( B(x)-P(x) \bigr)\, u, \\
        y &= \bigl( B(x)+P(x) \bigr)^\top z(x) + \bigl( S(x)-N(x) \bigr)\, u,
    \end{align}
\end{subequations}
associated with a Hamiltonian function $H\in\mathcal{C}^1(\mathcal{X},\mathbb{R})$, 
where $x(t)\in\mathcal{X}$ is the state, 
$u(t), y(t) \in \mathbb{R}^m$ are the input and output functions,
$r, z \in \mathcal{X}(\mathcal{S}, \mathbb{R]})$ are the \textit{time-flow and effort functions},
$E \in  \mathcal{C}(\mathcal{X}, \mathbb{R}^{n\times n})$ is the \textit{regular flow matrix}, 
$J,R \in  \mathcal{X}(\mathcal{S}, \mathbb{R}^{n\times n})$ are the \textit{structure and dissipation matrices},
$B, P \in \mathcal{C}(\mathcal{X},\mathbb{R}^{n \times p})$ are the \textit{port matrices} and $S, N \in  \mathcal{C}(\mathcal{X}, \mathbb{R}^{p \times p})$  are the \textit{feed-through matrices}.
Furthermore, the following property must hold:
\begin{equation*}
  \Gamma:= \begin{bmatrix}
    J & B \\ -B^\top & N
  \end{bmatrix} =- \Gamma^\top, \quad
  W:= \begin{bmatrix}
    R & P \\ P^\top & S
\end{bmatrix} \succeq 0,
\end{equation*}
i.e.\ $\Gamma$ is a skew-symmetric matrix and $W$ is a positive semidefinite matrix.
The gradient of the Hamiltonian function satisfies $\nabla H =E(x)^\top z$. The key property of such pHODE systems \eqref{pHODE} is the following \textit{dissipation inequality}
\begin{equation}\label{eq:dissineq}
\frac{d}{dt} H(x) = - \begin{pmatrix}
    z \\ u
\end{pmatrix}^\top W \begin{pmatrix}
    z \\ u
\end{pmatrix} + u^\top y,
\end{equation}
for the Hamiltonian (we refer to \cite[Section 2.2.1]{MehrmannMorandin} for a proof).

\begin{remark}
In many engineering applications, the feed-through matrices $S$ and $N$  may vanish. 
In this case, one has to also set the port matrix $P$ to zero because otherwise 
\begin{equation*}
   W:= \begin{bmatrix}
      R & P \\ P^\top & 0
\end{bmatrix} 
\end{equation*}
will become indefinite for $P \neq 0$. 
Therefore, we will often consider the special case of port-Hamiltonian ODE systems 
\begin{subequations}\label{pHODE.special}
    \begin{align}
        E(x) \dot{x} &= \bigl( J(x)-R(x) \bigr)\, z(x) + B(x)\, u, \\
        y &= B(x)^\top z(x).
    \end{align}
\end{subequations}
\end{remark}

\section{Coupling Port-Hamiltonian ODE Systems}\label{sec:coupling}
To define \textit{coupled port-Hamiltonian ODE systems}, we
consider $s$ subsystems of pHODEs~\eqref{pHODE.special} 
\begin{align*}
E_i(x_i)\, \dot{x}_i(t) + r_i(x_i)
    &= \bigl( J_i(x_i)  - R_i(x_i)\bigr) \,z_i(x_i)
         + \hat{B}_i(x_i)\,\hat{u}_i (t) + \bar{B}_i(x_i)\,\bar{u}_i (t), \\
   \hat{y}_i(t) &= \hat{B}_i(x_i)^\top z_i(x_i), \\   
   \bar{y}_i(t) &= \bar{B}_i(x_i)^\top z_i(x_i),
\end{align*}
$i=1,\ldots,s$, with $J_i=-J_i^\top$ and $R_i \succeq 0$ and a Hamiltonian $H_i(x_i)$, which fulfills 
the compatibility conditions $E_i^\top z_i = \nabla_{x_i} H_i$, $i=1,\ldots,s$. 
Here, we distinguish inputs $\hat{u}_i$ and outputs $\hat{y}_i$, which arise in the coupling between the subsystems, 
and inputs $\bar{u}_i$ and outputs $\bar{y}_i$, which couple the single system with the outer world.

The subsystems are coupled via \textit{external} inputs and outputs by
\begin{equation}
\label{eq:split-input-output} 
    \begin{pmatrix}
    \hat{u}_1 \\ \vdots \\ \hat{u}_k
    \end{pmatrix} 
    + \hat{C} 
    \begin{pmatrix}
    \hat{y}_1 \\ \vdots \\ \hat{y}_k
    \end{pmatrix} = 0, \qquad \hat{C} = - \hat{C}^\top. 
\end{equation}
These $s$ systems can be condensed into a large monolithic port-Hamiltonian descriptor system (pHDAE)
consisting of differential-algebraic equations of the type~\eqref{pHODE.special}
(cf.~\cite{Bartel23} for the linear case) 
\begin{subequations}\label{cond.ph.dae}
\begin{align}
   E(x) \dot{x} &= \bigl( J(x)-R(x) \bigr)\,z(x)
                        +  \bar{B}(x)\,\bar{u}, \\
    \bar{y} &=  \bar{B}(x)^\top z(x)),
\end{align}
\end{subequations}
with $J(x):=\diag(J_1(x_1),\ldots,J_s(x_s))- \hat{B}\hat{C} \hat{B}^\top$, 
$F=\diag\,(F_1,\ldots, F_s)$ 
for $F_i \in \{E, R, \bar{B}\}$ and $H(x):=\sum_{i=1}^s H_i(x_i)$.

The Schur complement $B C B^\top$ defines the off-diagonal block entries of $J$. 
\begin{remark}\label{remark1}
   If the matrix $\hat{C}$ in the input-output coupling~\eqref{eq:split-input-output} is not skew-symmetric, then one can split 
   $\hat{C} = \hat{C}^{\rm symm}+ \hat{C}^{\rm skew-symm}$ 
   into its symmetric and skew-symmetric parts, and one gets an equivalence with the pHODE system~\eqref{cond.ph.dae} by replacing $J(x)$ and $R(x)$ by 
   \begin{equation*}
      J(x) := \diag\bigl(J_1(x_1),\ldots,J_s(x_s)\bigr)- \hat{B} \hat{C}^{\rm skew-symm} \hat{B}^\top
   \quad
   \text{and} \quad  
      R(x):=\diag\bigl(R_1(x_1),\ldots,R_s(x_s)\bigr) + \hat{B} \hat{C}^{\rm symm} \hat{B}^\top,
   \end{equation*}
   provided that $R(x)\succeq 0$ is still positive semidefinite. 
   If not, one must use~\eqref{eq.equiv.dae} to follow the approach discussed in the following Remark~\ref{remark2}. 
\end{remark}
\begin{remark} \label{remark2}
In many engineering applications, the inputs and outputs of the subsystems are coupled via a skew-symmetric coupling condition. 
In general, if an arbitrary linear relation $M u + Ny =0$ exists, then the coupled system of pHODEs is equivalent to a pHDAE system (cf.\ \cite[Section 2.2.4]{MehrmannMorandin}), but this comes at the cost of introducing additional dummy variables $\hat{u}$ and $\hat{y}$ that copy $u$ and $y$:
\begin{subequations}\label{eq.equiv.dae}
\begin{align}
    \begin{bmatrix}
        E \dot{x} \\ 0 \\ 0 \\ 0
    \end{bmatrix}
    & =
    \begin{bmatrix} 
            \diag(J_1,\ldots, J_s)-\diag(R_1,\ldots, R_s) & 
            B & 0 & 0 \\ 
            -B^\top & 0 & I & -M^\top \\
            0 & - I & 0 & -N^\top \\
            0 & M & N & 0      
    \end{bmatrix}
    \begin{bmatrix}
        \nabla H(x) \\
        \hat{u} \\
        \hat{y} \\
        0
    \end{bmatrix}
    +
    \begin{bmatrix}
        0 \\ 0 \\ I \\ 0
    \end{bmatrix} u,\\
    y & = \hat{y}.
\end{align} 
\end{subequations}
\end{remark}

\section{Decoupling pHODE systems}\label{section-decoupling1}
Decoupling is only possible if there is a differentiable transformation of the state $w=F(x)$ such that the Hamiltonian becomes separable: 
\begin{equation*}
   H(x)=\tilde{H}(w)= \sum_{i=1}^s \tilde{H}_i(w_i), \quad
    w = \begin{bmatrix}
        w_1 \\ w_2 \\ \vdots \\ w_s
    \end{bmatrix}, \quad 
    n_i:=\dim(w_i),\:i=1,\ldots,s, \quad \sum_{i=1}^s n_i =n.
\end{equation*}
In this case, the pHODE system~\eqref{pHODE.special} can be transformed into the pHODE system
\begin{subequations}\label{ode.trafo}
\begin{align}
(F'(x)^{-1})^{\top} E(x) \dot x & = \underbrace{(F'(x)^{-1})^{\top}  E(x) F'(x)^{-1}}_{\displaystyle \tilde{E}(w)} \dot w \\
    & = \underbrace{ (F'(x)^{-1})^{\top} \bigl(J(x)-R(x)\bigr) F'(x)^{-1}}_{\displaystyle \tilde{J}(w) - \tilde{R}(w)} 
    \underbrace{F'(x) z(x)}_{\displaystyle \tilde{z}(w)} 
    + \underbrace{ (F'(x)^{-1})^{\top} B(x)}_{\displaystyle \tilde{B}(w)} u, 
    \nonumber\\
  y &= \tilde{B}(w)^\top \tilde{z}(w),
\end{align}
\end{subequations}
with $\tilde{E}^\top \tilde{z}(w) = \nabla_w \tilde{H}(w)$, as $\nabla_x H(x)= F'(x)^\top \nabla_w \tilde{H}(w) $ holds. 
Note that we have replaced $x$ with $F^{-1}(w)$.

Next, we split $\tilde{J}$ and $\tilde{R}$ with respect to splitting of $w$:
\begin{equation*}
\tilde{J} = \diag(\tilde{J}_{11},\ldots, \tilde{J}_{ss})  + \tilde{J}^{\rm offdiag}, \quad
\tilde{R} = \diag(\tilde{R}_{11},\ldots, \tilde{R}_{ss})\tilde{R}^{\rm diag} + \tilde{R}^{\rm offdiag}, 
\quad
\tilde{B} = \begin{bmatrix}
    \tilde{B}_1 \\ \vdots \\\tilde{B}_s
\end{bmatrix},
\end{equation*}
where we have used the abbreviations for the $n_i\times n_j$ matrices  $\tilde{J}_{ij}$ and  $\tilde{R}_{ij}$
\begin{equation*}
  \tilde{J}_{ij};=\left(J_{k,l}\right)_
  {\substack{k=\sum_{l=0}^{i-1} n_l +1,  \ldots , \sum_{l=0}^{i} n_l \\ l=\sum_{l=0}^{j-1} n_l +1,  \ldots , \sum_{l=0}^{j} n_l}} 
  \quad
   \tilde{R}_{ij}:=\left(R_{k,l}\right)_
  {\substack{k=\sum_{l=0}^{i-1} n_l +1,  \ldots , \sum_{l=0}^{i} n_l \\ l=\sum_{l=0}^{j-1} n_l +1,  \ldots , \sum_{l=0}^{j} n_l}}, \quad n_0:=0.
\end{equation*}
Provided that $\tilde{J}_{ii}, \tilde{R}_{ii}$ and $\tilde{B}_i$
depend only on $w_i$ (for $i=1,\ldots,s$), and the transformed flow matrix is block-diagonal of the form
$\tilde{E} = \diag(\tilde{E}_{11}(w_i),\ldots, \tilde{E}_{ss}(w_s))$, 
we can decouple the pHODE system into $s$ coupled pHODE systems as follows:

\subsection*{Case 1: $\tilde{R}^{\rm offdiag}=0$.} 
If we define $\hat{B}=I$ and $\hat{C}=-\tilde{J}^{\rm offdiag}$, then 
the pHODE system~\eqref{ode.trafo} is analytically equivalent to the $s$ pHODE subsystems 
\begin{subequations}\label{decoupled}
\begin{equation*} 
\begin{aligned}
    \tilde{E}_{ii}(w_i) \,\dot{w}_i 
    &= \bigl(\tilde{J}_{ii} -\tilde{R}_{ii}(w_i)\bigr) \,\tilde{z}_i(w_i) 
            + \begin{bmatrix}
              I , & \tilde{B}_i(w_i)
           \end{bmatrix} \begin{bmatrix}
               \hat{u}_i \\ u
           \end{bmatrix},
\qquad
  \begin{bmatrix}
      \hat{y}_i \\ \bar{y}_i
  \end{bmatrix} &= \begin{bmatrix}
              I, & \tilde{B}_i(w_i)
           \end{bmatrix}^\top \tilde{z}_i(w_i),
\end{aligned}
\end{equation*}
coupled by the skew-symmetric coupling condition 
\begin{equation*}
    \begin{pmatrix}
    \hat{u}_1 \\ \vdots \\ \hat{u}_k
    \end{pmatrix} 
    + \hat{C} 
    \begin{pmatrix}
    \hat{y}_1 \\ \vdots \\ \hat{y}_k
    \end{pmatrix} = 0, \qquad \hat{C} = - \hat{C}^\top,
\end{equation*}
and fulfilling the compatibility conditions $\tilde{E}_{ii}^\top \tilde{z}_i = \nabla_{w_i} \tilde{H}_{w_i}$.
\end{subequations}

\subsection*{Case 2: $\tilde{R}^{\rm offdiag} \neq0$.} 
Here we define $\hat{B}=I$ and $\hat{C}= -\tilde{J}^{\rm offdiag} + \tilde{R}^{\rm offdiag}$ to obtain the decoupled system~\eqref{decoupled}. 
However, the coupling condition is no longer skew-symmetric, but a general linear one with $M=I$ and $N=\hat{C}$ (see Remark~\ref{remark1}).

\subsection*{Case 3: General coupling.} 
We can use $s$ different port matrices $\hat{B}_i$ to define the equivalent pHODE system
\begin{subequations}\label{decoupled.case3}
\begin{equation} 
\begin{aligned}
   \tilde{E}_{ii} (w_i) \dot{w}_i 
    &= \bigl( \tilde{J}_{ii}(w_i)-\tilde{R}_{ii}(w_i) \bigr) \,\tilde{z}_i(w_i) 
            + \begin{bmatrix}
              \hat{B}_i, & \tilde{B}_i(w_i)
           \end{bmatrix} \begin{bmatrix}
               \hat{u}_i \\ u
           \end{bmatrix}, 
           \\
  \begin{bmatrix}
      \hat{y}_i \\ \bar{y}_i
  \end{bmatrix} &= \begin{bmatrix}
              \hat{B}_i , & \tilde{B}_i(w_i)
           \end{bmatrix}^\top \tilde{z}_i(w_i),
\end{aligned}
\end{equation}
coupled by the input-output coupling condition 
\begin{equation}
    \begin{pmatrix}
    \hat{u}_1 \\ \vdots \\ \hat{u}_s
    \end{pmatrix} 
    + \hat{C} 
    \begin{pmatrix}
    \hat{y}_1 \\ \vdots \\ \hat{y}_s
    \end{pmatrix} = 0, 
    \qquad
    \hat{C} := \begin{pmatrix}
    0 & \hat{C}_{12} & \cdots & \hat{C}_{1s} \\
    \hat{C}_{12}^\top & 0 & \cdots & \hat{C}_{2s} \\
    \vdots & \vdots & \ddots & \vdots \\
    \hat{C}_{1s}^\top & \hat{C}_{2s}^\top & \cdots & 0
\end{pmatrix}
\end{equation}
provided that a) 
\begin{equation*}
   \tilde{J}_{ij}-\tilde{R}_{ij} = -\hat{B}_i \hat{C}_{ij} \hat{B}_j^\top 
   \:\Rightarrow\: 
    \tilde{J}_{ij} = -\hat{B}_i \hat{C}^{\rm symm}_{ij} \hat{B}_j^\top 
    \:\text{and}\:  
    \tilde{R}_{ij} = \hat{B}_i \hat{C}^{\rm skew-symm}_{ij} \hat{B}_j^\top,
\end{equation*}
holds with $\hat{C}^{\rm symm}$ and $\hat{C}^{\rm skew-symm}$ corresponding to the symmetric and skew-symmetric part of $\hat{C}$, respectively, and
\end{subequations}
and b) the port matrices $\hat{B}_i$ depend only on $w_i$. 
The latter can always be achieved by defining $\hat{B}_i=I$ and 
\begin{equation*}
   \hat{C}_{ij} := -( \tilde{J}_{ij}-R_{ij}).
\end{equation*}
\begin{remark}
    In the case of a block-diagonal matrix $\tilde{R}$, i.e., $\tilde{R}_{ij} =0$ for $i \neq j$, we obtain a skew-symmetric coupling with
\begin{equation*}
\begin{pmatrix}
   0 & \hat{C}_{12} & \cdots & \hat{C}_{1p} \\
    -\hat{C}_{12}^\top & 0 & \cdots & \hat{C}_{2p} \\
   \vdots & \vdots & \ddots & \vdots \\
    -\hat{C}_{1p}^\top & -\hat{C}_{2p}^\top & \cdots & 0
\end{pmatrix}.    
\end{equation*}
  For this case, the obvious choice in the last section was 
  $\hat{B}_i=\hat{B}_j=I$, $\hat{C}_{ij}= -\tilde{J}_{ij}$.  
\end{remark}
\begin{remark}
 Since $R$ is positive semidefinite, all $R_{ii}$ are also positive semidefinite. 
 Additionally, if $J$ is skew-symmetric, then the block-diagonal matrices $J_{ii}$ are also skew-symmetric, 
 which defines $s$ pHODE systems. 
\end{remark}

\begin{example}
Consider the dissipative Hamiltonian system (pHODE system without external input, see Fig.~\ref{fig:two-masses_two-dampers_three-springs}) 
\begin{equation*}
    \dot{x} = (J-R) \,\nabla H(x),
\end{equation*}
    with
\begin{equation*}
    J = \begin{bmatrix}
        0 & -1 & -1 & 0 & 0 \\
        1 & 0 & 0 & 0 & 0 \\
        1 & 0 & 0 & -1 & 0 \\
        0 & 0 & 1 & 0 & -1 \\
        0 & 0 & 0 & 1 & 0 
    \end{bmatrix}, \quad
     R = \begin{bmatrix}
        r_1 & 0 & 0 & 0 & 0 \\
        0 & 0 & 0 & 0 & 0 \\
        0 & 0 & 0 & 0 & 0 \\
        0 & 0 & 0 & r_2 & 0 \\
        0 & 0 & 0 & 0 & 0 
    \end{bmatrix}, \quad
    H(x)= x^\top 
    \begin{bmatrix}
        \frac{1}{m_1} & 0 & 0 & 0 & 0 \\
        0 & K_1 & 0 & 0 & 0 \\
        0 & 0 & K & 0 & 0 \\
        0 & 0 & 0 & \frac{1}{m_2} & 0 \\
        0 & 0 & 0 & 0 & K_2 
    \end{bmatrix}
    x.
\end{equation*}
Next, we want to decouple this pHODE into two pHODE subsystems
      \begin{alignat*}{2}
        \dot{x}_1 & = (J_1-R_1) \,\nabla_1 H_1(x_1) + B_1 u_1, 
        & \qquad\ \qquad \qquad \quad 
                    \dot{x}_2 & = (J_2-R_2) \,\nabla_2 H_2(x_2) + B_2 u_2, \\
        y_1 & = B_1^\top \nabla_1 H_1(x_1), 
        & 
                    y_2 & = B_2^\top \nabla_2 H_2(x_2), \\
        &  & 0  = u+C y . \qquad \  &
    \end{alignat*}
    We have a lot of flexibility here: for example, if we split with respect to the first three and last two components, we have
    \begin{equation*}
       J_1= \begin{bmatrix}
           0 & - 1 & -1 \\
           1 & 0 & 0 \\
           1 & 0 & 0
       \end{bmatrix},
       \quad
       J_2=\begin{bmatrix}
           0 & -1\\
           1 & 0
       \end{bmatrix},
       \quad
       R_1=\begin{bmatrix}
           r_1 & 0 & 0 \\
           0 & 0 & 0 \\
           0 & 0 & 0
       \end{bmatrix},
       \quad
       R_2=\begin{bmatrix}
           r_2 & 0 \\
           0 & 0
       \end{bmatrix}
       \end{equation*}
       \begin{equation*}
       H_1(x_1) = x_1^\top \begin{bmatrix}
           \frac{1}{m_1} & 0 & 0\\
           0 & K_1 & 0 \\
           0 & 0 & K
       \end{bmatrix} x_1,
       \quad 
       H_2(x_2) = x_2^\top \begin{bmatrix}
           \frac{1}{m_1} & 0 \\
           0 & K_2
       \end{bmatrix} x_2.
       \end{equation*}
       The flexibility lies in the definition of the port matrices $B_1$ and $B_2$:
       \begin{itemize}
           \item[a)] We set $B_1=I_3$ and $B_2=I_2$ and
           \begin{equation*}
           C=-\begin{bmatrix}
        0 & 0 & 0 & 0 & 0 \\
        0 & 0 & 0 & 0 & 0 \\
        0 & 0 & 0 & -1 & 0 \\
        0 & 0 & 1 & 0 & 0 \\
        0 & 0 & 0 & 0 & 0 
    \end{bmatrix} = -J^{\rm offdiag}.
           \end{equation*}
           \item[b)] Here, we set 
           \begin{align*}
        &
           B_1 = \begin{bmatrix}
               0 \\ 0 \\ 1
          \end{bmatrix}, \quad 
        B_2
        = \begin{bmatrix}
               -1\\ 0 
        \end{bmatrix}, \quad
            C= \begin{bmatrix}
                0 & C_{12} \\ -C_{12}^\top & 0
            \end{bmatrix}, \quad
            C_{12} = -1,
        \end{align*}
    and check
    \begin{align*}
        B_1 C_{12} B_2^\top &=  \begin{bmatrix}
            0 \\ 0 \\ 1
        \end{bmatrix} \cdot (-1) \cdot
        \begin{bmatrix}
            -1 , & 0 
        \end{bmatrix} =
        \begin{bmatrix}
            0 & 0 \\
            0 & 0 \\
            1 & 0
        \end{bmatrix} = -J_{12} \quad \text{with}\;
        J= \begin{bmatrix}
            J_1 & J_{12} \\ -J_{12}^\top & J_2
        \end{bmatrix}.
    \end{align*}
    This splitting is in accordance with the modeling of a two-mass and three-spring system with damping, as described in~\cite{Bartel23}.
    This is achieved by setting $x_1 = [p_1, q_1, q_1 - q]^\top$ and $x_2 = [p_2, q_2]^\top$. 
\end{itemize}
    \begin{remark}
\noindent Properly choosing the port matrices minimizes the dimension of the coupling condition, i.e., the coupling matrix $\hat{C}$. 
With the former, we need only one scalar input for each subsystem. 
However, choosing $B_1$ and $B_2$ as identities requires as many inputs as there are state variables.

\end{remark}

\begin{figure}[H]
\begin{center}
\includegraphics[width=0.5\textwidth]{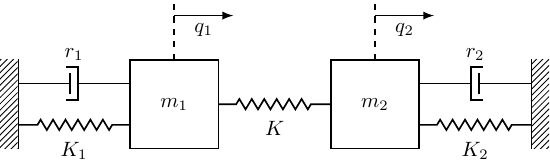}
\end{center}
\caption{\label{fig:two-masses_two-dampers_three-springs} ODE two masses oscillator with damping~\cite{Bartel23}. 
The coordinates $q_1,\, q_2$ describe the position of the masses.}
\end{figure}

       \noindent 
       Another splitting splits with respect to components 1,4 and components 2,3,5. In this case we have
        \begin{equation*}
       J_1=\begin{bmatrix}
           0 & 0 \\
           0 & 0
       \end{bmatrix}, \quad
       J_2=\begin{bmatrix}
           0 & 0 & 0 \\
           0 & 0 & 0 \\
           0 & 0 & 0
       \end{bmatrix}, \quad
       R_1=\diag(r_1,r_2), \quad
       R_2= \begin{bmatrix}
           0 & 0 & 0 \\
           0 & 0 & 0 \\
           0 & 0 & 0
       \end{bmatrix},
        \end{equation*}
         \begin{equation*}
       H_1(x_1)= x_1^\top \diag\Bigl( \frac{1}{m_1}, \frac{1}{m_2}\Bigr) x_1,\quad
       H_2=x_2^\top \diag(K_1,K,K_2)x_2.
        \end{equation*}
       A nontrivial choice of the port matrices is
    \begin{equation*}
       B_1=\begin{bmatrix}
           1 & 0 \\ 0 & -1
       \end{bmatrix},\quad
       B_2=\begin{bmatrix}
           -1 & 1 \\
           0 & 1 \\
           1 & 0
       \end{bmatrix},\quad
       C_{12}=\begin{bmatrix}
           0 & 1 \\ -1 & 0
       \end{bmatrix},
        \end{equation*}
       and we check $B_1 C_{12} B_2^\top = -J_{12}$. One may also split the system into 3, 4 and 5 subsystems as well.
\end{example}

\begin{example}[Poroelastic network model]
A discretized version of the poroelastic network model~\cite{Altmann21} reads
\begin{align*}
    E \dot{x} & = (J-R) x + B u, \\
    y & = B^\top x,
\end{align*}
with
\begin{equation*}
E=\begin{bmatrix}
    \rho M_u & 0 & 0 \\
    0 & K_u(\mu,\lambda) & 0 \\
    0 & 0 & \frac{1}{M} M_p
   \end{bmatrix}, \quad 
J= \begin{bmatrix}
    0 & - K_u(\mu,\lambda) & \alpha D^\top \\
    K_u(\mu,\lambda)^\top & 0 & 0 \\
    -\alpha D & 0 & 0
    \end{bmatrix}, \quad
R= \begin{bmatrix}
    0 & 0 & 0 \\
    0 & 0 & 0\\
    0 & 0 & \frac{\kappa}{\nu} K_p
\end{bmatrix}, 
\end{equation*}
\begin{equation*}
B= \begin{bmatrix}
    B_f & 0 \\ 0 & 0 \\ 0 & B_g
    \end{bmatrix},
    \quad
x = \begin{bmatrix}
    w(t) \\ v(t) \\ p(t)
     \end{bmatrix}, 
     \quad
u= \begin{bmatrix}
    f(t) \\ g(t)
    \end{bmatrix},
\end{equation*}
where $w$, $v$, $p$ are the discretized velocity, displacement and pressure fields, $f$, $g$ are the discretized volume-distributed forces and the external injection, $\mu$ and $\lambda$ are the Lame coefficients,
$\alpha$ is the Biot-Willes fluid-solid coupling coefficient, 
$\kappa$ is the permeability and $\rho$ is the density. 
$M_u$ and $M_p$ are mass matrices with respect to \ $u$ and $p$, $K_u$ is a stiffness matrix with respect to \ $u$ and $D$ is a damping matrix, see \cite{Benchmark} for more details. We can easily split this implicit pHODE system with respect to \ the first two and the third block and obtain
\begin{subequations}
    \begin{align}
    \begin{bmatrix}
    \rho M_u & 0  \\
    0 & K_u(\mu,\lambda) 
\end{bmatrix} \cdot
\begin{bmatrix}
    \dot{w}(t) \\ \dot{u}(t) 
\end{bmatrix} & = \begin{bmatrix}
    0 & - K_u(\mu,\lambda)  \\
    K_u(\mu,\lambda)^\top & 0 
\end{bmatrix} \cdot
\begin{bmatrix}
    w(t) \\ v(t) 
\end{bmatrix}
+ \begin{bmatrix}
    \alpha D^\top \\ 0
\end{bmatrix} u_1, \\
y_1 & = \begin{bmatrix}
    \alpha D^\top \\ 0
\end{bmatrix}^\top \begin{bmatrix}
    w(t) \\ v(t) 
\end{bmatrix} =  \alpha D w(t),
\end{align}
\begin{align}
\frac{1}{M} M_p \dot{p} & = 
-\frac{\kappa}{\nu} K_p p + u_2, \\
y_2 & = p,
    \end{align}
    \begin{align} 0 & = 
        \begin{bmatrix}
            u_1 \\ u_2
        \end{bmatrix} +
        \begin{bmatrix}
            0 & I \\ -I & 0
        \end{bmatrix}
        \begin{bmatrix}
            y_1 \\ y_2
        \end{bmatrix}.
    \end{align}
\end{subequations}
Because the first subsystem lacks a dissipative part and the second lacks a skew-symmetric part, symplectic time integration can be used for the first part and dissipative time integration can be used for the second part in a dynamic iteration procedure. 
\end{example}
\begin{example}[Discretized full set of Maxwell's equations]
The full set of Maxwell's equations formulated with electrodynamic potentials and discretized forms a pHDAE system with 
\begin{equation} \label{eq:FullMaxwell-pHDAE}
   \underbrace{
   \begin{bmatrix}
   \fMeps                 &            \fMeps \Gr& 0\\
   \Gr^\top \fMeps      & \Gr^\top \fMeps \Gr& 0\\
   0 & 0  & \fMmu
   \end{bmatrix}
   }_{:=\fitmat{E}_{\rm Maxwell} } 
   d_t
\underbrace{
   \begin{bmatrix}
     d_t\fitvec{a}         \\
     \fitvec{\varphi}   \\ 
     \fitvec{h} 
   \end{bmatrix}
   }_{:=\fitvec{x}}
   =
   \left[
   \smash[b]{\underbrace{ 
   \begin{bmatrix}
   0      & 0 & -\C^\top\\
   0    & 0 &       0\\
   \C & 0  & 0
   \end{bmatrix}
   }_{:= \fitmat{J} (= -\fitmat{J}^\top)}
   }
   -
   \smash[b]{
   \underbrace{
   \begin{bmatrix}
           \fMkap &         \fMkap \Gr  & 0\\
   \Gr^\top \fMkap & \Gr^\top\fMkap \Gr  & 0\\
   0    & 0 &       0
   \end{bmatrix}
   }_{=R \ge 0}
   }
   \right]
   \underbrace{
   \begin{bmatrix}
    d_t\fitvec{a}  \\
    \fitvec{\varphi} \\ 
    \fitvec{h} 
   \end{bmatrix}
   }_{:=\fitvec{z}}
   +
   \begin{bmatrix}
   \fitmat{I} \\
   \Gr^\top \\
    0         
   \end{bmatrix}
   \fitvec{j}_{\mathrm{s}},
\end{equation}
where $\C$, $\Gr$ are discrete curl and gradient matrices, $\fMkap$, $\fMeps$, $\fMmu$ are discrete material matrices for electric conductivity $\kappa$, permittivity $\varepsilon$ and permeability $\mu$, respectively, $\fitvec{a}$, $\fitvec{\varphi}$,
$\fitvec{h}$ are the degree of freedom vectors and $\fitvec{\bm u} = \fitvec{j}_{\mathrm{s}}$ is the input vector of source currents and the output vector $\fitvec{y}= d_t \fitvec{a}+\Gr{\varphi} (=- \fitvec{e})$ is the negated vector of electric grid voltages. 
The Hamiltonian associated with \eqref{eq:FullMaxwell-pHDAE} is the discrete electromagnetic grid energy is given by
\begin{equation}\label{eq:Hamiltonian_Maxwell}
    H= \frac{1}{2}
    \left(\fitvec{\varphi}^\top \Gr^\top \fMeps\Gr\fitvec{\varphi}
    +  \fitvec{h}^\top \fMmu\fitvec{h}\right) \\
    + {{\varphi}^\top \Gr^\top \fMeps {\it d_t} \fitvec{a}} 
    + d_t\fitvec{a}^\top \fMeps d_t\fitvec{a},
\end{equation}
 and the compatibility condition holds with 
\begin{equation}
\nabla H = 
\left(
\fMeps(d_t\fitvec{a} + \Gr {\varphi}), \Gr^\top \fMeps ({\it d_t}\fitvec{a} + \Gr {\varphi} 
),\fMmu\fitvec{h}
\right)^\top 
= \fitmat{E}_{\mathrm{Maxwell} }^\top \fitvec{x}.  
\end{equation}
Again, we can easily split with respect to \ the first two and the last block with Hamiltonians 
\begin{equation*}
H_1= \frac{1}{2}
    \fitvec{\varphi}^\top \Gr^\top \fMeps\Gr\fitvec{\varphi}
    + {{\varphi}^\top \Gr^\top \fMeps {\it d_t} \fitvec{a}} 
    + d_t\fitvec{a}^\top \fMeps d_t\fitvec{a}\quad\text{and}\quad
    H_2= \frac{1}{2}
      \fitvec{h}^\top \fMmu\fitvec{h},
\end{equation*}
\[
  \begin{bmatrix}
   \fMeps                 &            \fMeps \Gr\\
   \Gr^\top \fMeps      & \Gr^\top \fMeps \Gr
   \end{bmatrix}
   d_t
   \begin{bmatrix}
    d_t\fitvec{a}       \\
    \fitvec{\varphi} 
   \end{bmatrix}
   =
     \begin{bmatrix}
           \fMkap &         \fMkap \Gr  \\
   \Gr^\top\fMkap & \Gr^\top\fMkap \Gr  
   \end{bmatrix}
   \begin{bmatrix}
   d_t\fitvec{a}   \\
   \fitvec{\varphi}  
   \end{bmatrix}
   +
   \begin{bmatrix}
   \fitmat{I}  \\
   \Gr^\top
   \end{bmatrix}
    \fitvec{j}_{\mathrm{s}} + 
   \begin{bmatrix}
      \fitmat{C}^\top \\
      \fitmat{0} 
   \end{bmatrix}
   \fitmat{u_1},
\]
\[
   \fitmat{\tilde{y}}_1  =  
   \begin{bmatrix}
     \fitmat{I}  \\
    \fitmat{G}^\top 
   \end{bmatrix}^\top 
   \begin{bmatrix}
    d_t\fitvec{a}  \\
   \fitvec{\varphi}  
   \end{bmatrix}, \quad 
   \fitmat{y}_1  = 
   \begin{bmatrix}
     \fitmat{C}^\top \\
    \fitmat{0}
   \end{bmatrix}^\top 
   \begin{bmatrix}{c}
   d_t\fitvec{a}    \\
    \fitvec{\varphi}  
   \end{bmatrix},
\]
\[
   \fMmu d_t
\fitvec{h}
 =  \fitvec{u}_2, \quad \fitvec{y}_2  = \fitvec{h},
\quad 
  0 =
        \begin{bmatrix}
            \fitvec{u}_1 \\ \fitvec{u}_2
        \end{bmatrix} +
        \begin{bmatrix}
            \fitmat{0} & \fitmat{I} \\ -\fitmat{I} & \fitmat{0}
        \end{bmatrix}
        \begin{bmatrix}
            \fitvec{y}_1 \\ \fitvec{y}_2
        \end{bmatrix}.
\]
\end{example}
\section{Numerical Aspects of Monolithic and Decoupled Formulations}\label{sec:advantages}

In numerical simulations, the monolithic and decoupled formulations of port-Hamiltonian ODE systems offer distinct advantages.
The monolithic formulation preserves the full system structure and typically exhibits superior stability properties, particularly when integrated with structure-preserving or symplectic time-stepping schemes. 
This makes it ideal for stiff or highly coupled systems where global energy conservation is crucial. 
Nevertheless, one may reduce computational time by applying operator splitting techniques to the general systems which preserve the port-Hamiltonian structure, for example, $J-R$ coupling~\cite{Bartel25} or diagonal/off-diagonal coupling~\cite{Lorenz2025}, jointly with multirate approaches.

In contrast, the decoupled formulation offers greater flexibility and computational efficiency in multiphysics or large-scale contexts. 
It allows for the use of specialized numerical solvers, such as mimetic or energy-consistent methods, which are tailored to individual subsystems and respect local conservation laws. 
This partitioned approach also enables parallel or distributed computation via dynamic iteration~\cite{Guenther22}, although with potentially weaker global stability.
The monolithic formulation favors global stability, whereas the decoupled approach provides modularity and solver flexibility, at the cost of potentially reduced overall stability.

\section{Conclusion}\label{sec:conclusion}
We have presented a unified framework for coupling and decoupling nonlinear port-Hamiltonian ordinary differential equation systems. 
The analysis demonstrates how multiple pH subsystems can be reformulated as a single monolithic system, and conversely, how a monolithic pHODE can be decomposed into weakly coupled subsystems without losing the underlying Hamiltonian structure.
The results provide a theoretical basis for structure-preserving simulation strategies that balance stability and flexibility.
Future research will address adaptive coupling strategies, structure-preserving time discretizations, and extensions to port-Hamiltonian differential-algebraic \cite{Guenther21} and stochastic systems \cite{DiPersio25}.

\section*{Acknowledgements}
M.~Ehrhardt and M.~Günther acknowledge funding by the Deutsche Forschungsgemeinschaft (DFG, German Research Foundation) – Project-ID 531152215 – CRC 1701. 
D~\v{S}ev\v{c}ovi\v{c} was supported by Slovak research grant VEGA 1-0493-24.

\medskip
\noindent{\bf Data availability}
\\
No data was used for the research described in the article.

\end{document}